\documentclass{amsart}

\newtheorem{theorem}{Theorem}
\newtheorem{lemma}{Lemma}

\newtheorem{cor}{Corollary}
\newtheorem{prop}{Proposition}

\usepackage{amsmath}
\usepackage[psamsfonts]{amssymb}
\usepackage[mathscr]{euscript}

\newcommand{\cal}{\mathcal}

\title[Harmonic HS on Riemannian manifolds with skew-torsion \hfill]
{Harmonic  Hermitian structures on Riemannian manifolds with
skew-torsion }

\author{Johann Davidov}
\thanks{The author is partially supported by the National Science Fund, Ministry of Education and Science of Bulgaria under contract DN 12/2.}

\address{Johann Davidov\\Institute of Mathematics and Informatics \\
Bulgarian Academy of Sciences\\ Acad. G. Bonchev Str. Bl. 8\\ 1113 Sofia\\
Bulgaria.}

\begin{document}

\begin{abstract}

We find geometric conditions on a four-dimensional  Hermitian
manifold endowed with a metric connection with totally
skew-symmetric torsion under which the  complex structure is a
harmonic map from the manifold into its twistor space considered
with a natural family of Riemannian metrics defined by means of the
metric and the given connection on the base manifold.

\vspace{0,1cm} \noindent 2020 {\it Mathematics Subject
Classification.} 53C43, 58E20.

\vspace{0,1cm} \noindent {\it Keywords: Almost complex structures,
harmonic maps, twistor spaces, skew-symmetric torsion.}

\end{abstract}

\thispagestyle{empty}

\maketitle

\section{Intorduction}

\medskip

Recall that an almost complex structure on a Riemannian manifold is
called almost Hermitian or compatible if it is an orthogonal
endomorphism of the tangent bundle of the manifold. It is well known
that if a Riemannian manifold $(N,h)$ admits an  almost Hermitian
structure, it possesses many such structures. In \cite{D17, DHM},
for example, this is shown by considering an almost Hermitian
structure on $(N,h)$ as a section of the twistor bundle
$\pi:{\mathcal Z}\to N$ whose fibre at a point $p\in N$ consists of
all $g$-orthogonal complex structures $I_p:T_pN\to T_pN$
($I_p^2=-Id$) on the tangent space of $N$ at $p$. Thus, it is
natural to look for \textquotedblleft reasonable" criteria that
distinguish some of the almost Hermitian structures on a Riemannian
manifold. Since we shall consider the almost Hermitian structures on
$(N,h)$ as sections of the bundle ${\mathcal Z}$, let us note that
the fibre of this bundle is the compact Hermitian symmetric space
$O(2m)/U(m)$, $2m=dim\,M$, and its standard metric
$G=-\frac{1}{2}Trace\,I_1\circ I_2$ is K\"ahler-Einstein. The
Levi-Civita connection of $(N,h)$ gives rise to a splitting
$T{\mathcal Z}={\mathcal H}\oplus {\mathcal V}$ of the tangent
bundle of ${\mathcal Z}$ into horizontal and vertical parts. This
decomposition allows one to define a $1$-parameter family of
Riemanian metrics $h_t=\pi^{\ast}h+tG$, $t>0$, for which the
projection map $\pi:({\mathcal Z},h_t)\to (N,h)$ is a Riemannian
submersion with totally geodesic fibres. In the terminology of
\cite[Definition 9.7]{Besse}, the family $h_t$ is the canonical
variation of the metric $h$. Motivated by harmonic map theory, C. M.
Wood \cite{W1,W2} has suggested to consider as \textquotedblleft
optimal" those almost Hermitian structures $J:(N,h)\to ({\mathcal Z
},h_1)$ that are critical points of the energy functional under
variations through sections of ${\mathcal Z}$. In general, these
critical points are not harmonic maps, but, by analogy, they are
referred to as \textquotedblleft harmonic almost complex structures"
in \cite{W2}; they are also called \textquotedblleft harmonic
sections" in \cite{W1}, a term which seems more appropriate in the
context of this paper. Forgetting the bundle structure of ${\mathcal
Z}$, we can consider the almost Hermitian structures that are
critical points of the energy functional under variations through
all maps $N\to{\mathcal Z}$. These structures are genuine harmonic
maps from $(N,h)$ into $({\mathcal Z}, h_t)$, and we refer to
\cite{EL} for basic facts about such maps. This point of view is
taken in \cite{DHM} where the problem of when an almost Hermitian
structure on a Riemannian four-manifold is a harmonic map from the
manifold into its twistor space is discussed.

Every metric connection $D$ on $(N,h)$ yields a decomposition of
$T{\mathcal Z}$ into horizontal and vertical subbundles,  and having
such a connection we can define the corresponding family $h_t$ of
Riemannian metrics on the twistor space ${\mathcal Z}$. Of special
interest are metric connections with (totally) skew-symmetric
torsion which have many applications in various areas of mathematics
and physics. From the point of view of twistor theory, it is worth
noting that the Atiyah-Hitchin-Singer \cite{AHS} almost complex
structures on the twistor space of an oriented four-dimensional
Riemannian manifold defined by means of the Levi-Civita connection
and by a metric connection with skew-symmetric torsion have the same
integrability condition, see, for example, \cite{D20}. Now, recall
that a smooth map $f:(N,h)\to (N',h')$ of Riemannian manifolds is
harmonic if and only if the trace of its second fundamental form
$II_f(\nabla,\nabla')$ defined by means of the Levi-Civita
connections $\nabla$ of $(N,h)$ and $\nabla'$ of $(N',h')$ vanishes.
If $D$ and $D'$  are metric connections with skew-symmetric torsion
 a simple calculation shows that the trace of the second
fundamental form $II_f(D,D')$ defined by means of $D$ and $D'$
coincides with $Trace\,II_f(\nabla,\nabla')$. Thus, the metric
connections with skew-symmetric torsion can be used for studying
harmonic maps between Riemannian manifolds endowed with additional
structures preserving by such connections.

 In the present paper we
consider the twistor space of a four-dimensional Riemannian manifold
$(M,g)$ endowed with a metric connection $D$ with skew-symmetric
torsion and, by means of this connection, we define  the metrics
$h_t$. We find the conditions under which a Hermitian structure $J$
on $(M,g)$ is a harmonic map $J:(M,g)\to ({\mathcal Z},h_t)$. In
particular, if $D$ is the Bismut-Strominger connection
\cite{Bi,Str},  the unique metric connection with skew-torsion
preserving the complex structure, then  the map $J$ is always
harmonic. Two examples illustrating the obtained result are given,
one of them showing that, for a fixed metric, there can be many
metric connections with skew-torsion for which $J$ is a harmonic
map.

\section{Harmonic maps between Riemannian manifolds endowed with
metric connection with skew-torsion}

Let $(N,h)$, $(N',h')$ be Riemannian manifolds and $f:N\to N'$ a
smooth map. If $f^{\ast}TN'$ is the pull-back bundle of the bundle
$TN'$ under the map $f$, we can consider the differential
$f_{\ast}:TN\to TN'$ as a section of the bundle
$Hom(TN,f^{\ast}TN')$. Suppose we are given connections $D$ and $D'$
on $TN$ and $TN'$, respectively. Denote by $D'^{\ast}$ the
connection on $f^{\ast}TN'$ induced by the connection $D'$ of $TN'$.
The connections $D$ and $D'^{\ast}$ give rise to a connection
$\widehat D$ on the bundle $Hom(TN,f^{\ast}TN')$.  Define a bilinear
form on $N$ setting
$$
II_{f}(D,D')(X,Y)=(\widehat D_{X}f_{\ast})(Y),\quad X,Y\in TN.
$$
If $D$ and $D'$ are torsion-free, this form is symmetric. Recall
that the map $f$ is said to be harmonic if for the Levi-Civita
connections $\nabla$ and $\nabla'$ of $TN$ and $TN'$,
$$
Trace_{h} II_{f}(\nabla,\nabla')=0.
$$
Now, suppose that $D$ and $D'$ are metric connections with (totally)
skew-symmetric torsions $T$ and $T'$, i.e., the trilinear form
$$
{\mathcal T}(X,Y,Z)=h(T(X,Y),Z), \quad X,Y,Z\in TM,
$$
is skew-symmetric, and similarly for ${\mathcal
T}'(X,Y,Z)=h'(T'(X,Y),Z)$. Recall that, on a Riemannian manifold,
there is a unique metric connection with a given torsion; for an
explicit formula see, for example, \cite[Sec. 3.5, formula
(14)]{GKM}. Since the torsion $3$-forms ${\mathcal T}$ and
${\mathcal T}'$ are skew-symmetric,
\begin{equation}\label{D-nabla}
D_{X}Y=\nabla_{X}Y+\frac{1}{2}T(X,Y),\quad
D'_{X}Y=\nabla'_{X}Y+\frac{1}{2}T'(X,Y).
\end{equation}

\begin{lemma}\label{DLC-D}
Let $f:(N,h)\to (N',h')$ be a smooth map of Riemannian manifolds
endowed with metric connections $D$ and $D'$ with skew-symmetric
torsions $T$ and $T'$. Denote the Levi-Civita connections of $(N,h)$
and $(N,h')$ by $\nabla$ and $\nabla'$. Then
$$
Trace_{h}II_{f}(D,D')=Trace_{h}II_{f}(\nabla,\nabla').
$$
\end{lemma}

\noindent {\bf Proof}. Take a point $p\in N$ and let $E_1,...,E_n$
and $E_1',...,E_{n'}'$ be frames of $TN$ and $TN'$ defined on
neighbourhoods $U$ of $p$ and $U'$ of $p'=f(p)$ such that
$f(U)\subset U'$. Choose the frame $E_1,...,E_n$ to be
$h$-orthonormal. Then $f_{\ast}\circ E_i=\sum\limits_{\alpha=1}^{n'}
\lambda_{i\alpha}E'_{\alpha}\circ f$ where $\lambda_{i\alpha}$ are
smooth functions on $U$. Under this notation,
$$
\begin{array}{c}
Trace_{h}II_{f}(D,D')=\sum\limits_{i}[D'^{\ast}_{E_i}(f_{\ast}\circ
E_i)-f_{\ast}(D_{E_i}E_i)]\\[8pt]
=\sum\limits_{i}\sum\limits_{\alpha}[E_i(\lambda_{i\alpha})E'_{\alpha}\circ
f+\sum\limits_{\beta}\lambda_{i\alpha}\lambda_{i\beta}(D'_{E'_{\beta}}E'_{\alpha})\circ
f]-\sum\limits_{i}f_{\ast}(D_{E_i}E_i)\\[8pt]
=\sum\limits_{i}\sum\limits_{\alpha}[E_i(\lambda_{i\alpha})E'_{\alpha}\circ
f+\sum\limits_{\beta}\lambda_{i\alpha}\lambda_{i\beta}(\nabla'_{E'_{\beta}}E'_{\alpha})\circ
f]+\displaystyle{\frac{1}{2}}\sum\limits_{i}T'(f_{\ast}(E_i),f_{\ast}(E_i))\\[8pt]
-\sum\limits_{i}f_{\ast}(\nabla_{E_i}E_i)-\displaystyle{\frac{1}{2}}\sum\limits_{i}T(E_i,E_i)\\[8pt]
=Trace_{h}II_{f}(\nabla,\nabla').
\end{array}
$$

%\hfill $\Box$

\begin{cor}\label{har-D}
A map $f:N\to N'$ is harmonic if and only if
$$Trace_{h}II_{f}(D,D')=0.$$
\end{cor}

\section{Basics on twistor spaces}

We recall first some basic facts about twistor space following \cite
{D17,D20}.

\smallskip

Let $(M,g)$ be an oriented  (connected) Riemannian manifold of
dimension four. The metric $g$ induces a metric on the bundle of
two-vectors $\pi:\Lambda^2TM\to M$ by the formula
\begin{equation}\label{g on Lambda2}
g(v_1\wedge v_2,v_3\wedge v_4)=\frac{1}{2}det[g(v_i,v_j)]
\end{equation}
(the choice of the factor $1/2$ is explained below).

Let $\ast:\Lambda^kTM\to \Lambda^{4-k}M$, $k=0,...,4$, be the Hodge
star operator. Its restriction to $\Lambda^2TM$ is an involution,
thus we have the orthogonal decomposition
$$
\Lambda^2TM=\Lambda^2_{-}TM\oplus\Lambda^2_{+}TM,
$$
where $\Lambda^2_{\pm}TM$ are the subbundles of $\Lambda^2TM$
corresponding to the $(\pm 1)$-eigenvalues of the operator $\ast$.

\smallskip

Let $(E_1,E_2,E_3,E_4)$ be a local oriented orthonormal frame of
$TM$. Set
\begin{equation}\label{s-basis}
s_1^{\pm}=E_1\wedge E_2\pm E_3\wedge E_4, \quad s_2^{\pm}=E_1\wedge
E_3\pm E_4\wedge E_2, \quad s_3^{\pm}=E_1\wedge E_4\pm E_2\wedge
E_3.
\end{equation}
Then $(s_1^{\pm},s_2^{\pm},s_3^{\pm})$ is a local orthonormal frame
of $\Lambda^2_{\pm}TM$. This frame defines an orientation on
$\Lambda^2_{\pm}TM$ which does not depend on the choice of the frame
$(E_1,E_2,E_3,E_4)$ (see, for example, \cite{D17}). We call this
orientation "canonical".

For every $a\in\Lambda ^2TM$, define a skew-symmetric endomorphism
of $T_{\pi(a)}M$ by
\begin{equation}\label{cs}
g(K_{a}X,Y)=2g(a, X\wedge Y), \quad X,Y\in T_{\pi(a)}M.
\end{equation}
Denoting by $G$ the standard metric $-\frac{1}{2}Trace\,PQ$ on the
space of skew-symmetric endomorphisms, we have $G(K_a,K_b)=2g(a,b)$
for $a,b\in \Lambda ^2TM$. If $a\in\Lambda^2TM$ is of unit length,
then $K_a$ is a complex structure on the vector space $T_{\pi(a)}M$
compatible with the metric $g$, i.e., $g$-orthogonal. Conversely,
the $2$-vector $a$ dual to one half of the fundamental $2$-form of
such a complex structure is a unit vector in $\Lambda^2TM$.
Therefore the unit sphere bundle ${\mathcal Z}$ of $\Lambda ^2TM$
parametrizes the complex structures on the tangent spaces of $M$
compatible with the metric $g$ (so, the factor $1/2$ in (\ref{g on
Lambda2}) is chosen in order to have spheres with radius $1$). This
bundle is called the twistor space of the Riemannian manifold
$(M,g)$. Since $M$ is oriented, the manifold ${\mathcal Z}$ has two
connected components ${\mathcal Z}_{\pm}$ called the positive and
the negative twistor spaces of $(M,g)$. These are the unit sphere
subbundles of $\Lambda^2_{\pm}TM$. The bundle ${\mathcal Z}_{\pm}\to
M$ parametrizes  the complex structures on the tangent spaces of $M$
compatible with the metric and  $\pm$ the orientation via the
correspondence ${\mathcal Z}_{\pm}\ni\sigma\to K_{\sigma}$. Note
that changing the orientation of $M$ interchanges the roles of
$\Lambda^2_{-}TM$ and $\Lambda^2_{+}TM$ and, respectively, of
${\mathcal Z}_{-}$ and ${\mathcal Z}_{+}$.

\smallskip

The vertical space ${\mathcal V}_{\sigma}=\{V\in T_{\sigma}{\mathcal
Z}_{\pm}:~ \pi_{\ast}V=0\}$ of the bundle $\pi:{\mathcal Z}_{\pm}\to
M$ at a point $\sigma$ is the tangent space to the fibre of
${\mathcal Z}_{\pm}$ through $\sigma$. Thus, considering
$T_{\sigma}{\mathcal Z}_{\pm}$ as a subspace of
$T_{\sigma}(\Lambda^2_{\pm}TM)$, ${\mathcal V}_{\sigma}$ is the
orthogonal complement of ${\Bbb R}\sigma$ in
$\Lambda^2_{\pm}T_{\pi(\sigma)}M$.

%If $V\in{\mathcal V}_{\sigma}$ the skew-symmetric endomorphisms
%$K_{V}$ of $T_{\pi(\sigma)}M$ anti-commutes with $K_{\sigma}$.

Let $D$ be a metric connection on $(M,g)$. The induced connection on
$\Lambda^2TM$ will also be denoted by $D$. If
$(s_1^{\pm},s_2^{\pm},s_3^{\pm})$ is the orthonormall frame of
$\Lambda ^2_{\pm}TM$ defined by means of an oriented orthonormal
frame $(E_1,...,E_4)$ of $TM$ via (\ref{s-basis}), we have
$g(D_{X}s_i^{+},s_j^{-})= g(D_{X}s_i^{-},s_j^{+})=0$ and
$g(D_{X}s_i^{\pm},s_j^{\pm})=-g(D_{X}s_j^{\pm},s_i^{\pm})$ for every
$X\in TM$ and every $i,j=1,2,3$. Therefore the connection $D$
preserves the bundles $\Lambda^2_{\pm}TM$ and induces a metric
connection on each of these bundles denoted again by $D$. Let
$\sigma\in{\mathcal Z}_{\pm}$,  and let $s$ be a local section of
${\mathcal Z}_{\pm}$ such that $s(p)=\sigma$ where $p=\pi(\sigma)$.
Considering $s$ as a section of $\Lambda^2_{\pm}TM$, we have
$D_{X}s\perp s(p)$, i.e., $D_{X}s\in{\cal V}_{\sigma}$ for every
$X\in T_pM$ since $s$ has a constant length. Moreover,
$X^h_{\sigma}=s_{\ast}X-D_{X}s\in T_{\sigma}{\mathcal Z}_{\pm}$ is
the horizontal lift of $X$ at ${\sigma}$ with respect the connection
$D$ on $\Lambda^2_{\pm}TM$. Thus, the horizontal distribution of
$\Lambda^2_{\pm}TM$ with respect to $D$ is tangent to the twistor
space ${\mathcal Z}_{\pm}$. In this way, we have the decomposition
$T{\mathcal Z}_{\pm}={\mathcal H}\oplus {\mathcal V}$ of the tangent
bundle of ${\mathcal Z}_{\pm}$ into horizontal and vertical
components. The horizontal and vertical parts of a tangent vector
$A\in T{\mathcal Z}_{\pm}$ will be denoted by ${\mathcal H}A$ and
${\mathcal V}A$, respectively.

\smallskip

\noindent {\bf Convention}. In what follow, we shall freely identify
the bundle $\Lambda ^2TM$ with the bundle $A(TM)$ of
$g$-skew-symmetric endomorphism of $TM$ by means of the isomorphism
$a\to K_a$ defined by (\ref{cs}).

\smallskip

Using the basis (\ref{s-basis}), it is easy to check that if
$a,b\in\Lambda^2_{\pm}T_pM$, the isomorphism $\Lambda^2TM\cong
A(TM)$ sends $a\times b$ to $\pm\frac{1}{2}[K_a,K_b]$. In the case
when $a\in\Lambda^2_{+}T_pM$, $b\in\Lambda^2_{-}T_pM$, the
endomorphisms $K_a$ and $K_b$ of $T_pM$ commute. If
$a,b\in\Lambda_{\pm}T_pM$,
$$
K_{a}\circ K_{b}=-g(a,b)Id\pm K_{a\times b}.
$$
In particular, $K_a$ and $K_b$ anti-commute if and only if $a$ and
$b$ are orthogonal.

\smallskip

For every $t>0$, define a Riemannian metric $h_t=h_t^D$ on
${\mathcal Z}_{+}$ by
$$
h_t(X^h_{\sigma}+V,Y^h_{\sigma}+W)=g(X,Y)+tg(V,W)
$$
for $\sigma\in{\mathcal Z}_{+}$, $X,Y\in T_{\pi(\sigma)}M$,
$V,W\in{\mathcal V}_{\sigma}$.

\smallskip

\noindent {\bf Notation}. We set ${\mathcal Z}={\mathcal Z}_{+}$.
The sections $s_i^{+}$ of $\Lambda^2_{+}TM$ defined via
(\ref{s-basis}) will be denoted by $s_i$, $i=1,2,3$. The Levi-Civita
connection of the metric $h_t$ on ${\mathcal Z}$ will be denote by
$\widetilde{D}$.

\smallskip

By the Vilms theorem (see, for example, \cite[Theorem 9.59]{Besse}),
the projection map $\pi:({\mathcal Z},h_t)\to (M,g)$ is a Riemannian
submersion with totally geodesic fibres.  This can also be proved by
an easy  direct computation as is shown below.

\smallskip
\noindent {\bf Convention}. For the curvature tensor we adopt the
following definition $R^D(X,Y)=D_{[X,Y]}-[D_X,D_Y]$. The curvature
tensor of the induced connection on $\Lambda^2_{\pm}TM$ will again
be denoted by $R^D$. The curvature operator ${\mathcal
R}^D:\Lambda^2TM\to\Lambda^2TM$ is defined by $g({\mathcal
R}^D(X\wedge Y),Z\wedge U)=g(R^D(X,Y)Z,U)$.
\smallskip

The following formula is easy to be checked.
\begin{lemma}\label{CurOp}
If $a,b\in\Lambda^2_{+}T_pM$, then
$$
g(R^D(X,Y)a,b)=g({\mathcal R}^D(X\wedge Y),a\times b)
$$
for every $X,Y\in T_pM$.
\end{lemma}

Let $({\mathscr N},x_1,...,x_4)$ be a local coordinate system of $M$
and let $(E_1,...,E_4)$ be an oriented orthonormal frame of $TM$ on
${\mathscr N}$. If $(s_1,s_2,s_3)$ is the local frame of
$\Lambda^2_{+}TM$ define by (\ref{s-basis}), then $\widetilde
x_{\alpha}=x_{\alpha}\circ\pi$, $y_j(\sigma)=g(\tau,
(s_j\circ\pi)(\sigma))$, $1\leq \alpha \leq 4$, $1\leq j\leq 3$, are
local coordinates of $\Lambda^2_{+}TM$ on $\pi^{-1}({\mathscr N})$.

   The horizontal lift $X^h$ on $\pi^{-1}({\mathscr N})$ of a vector field
$$
X=\sum_{a=1}^4 X^{\alpha}\frac{\partial}{\partial x_{\alpha}}
$$
is given by
\begin{equation}\label{hl}
X^h=\sum_{\alpha=1}^4 (X^{\alpha}\circ\pi)\frac{\partial}{\partial
\widetilde{x}_{\alpha}}
-\sum_{j,k=1}^3y_j(g(D_{X}s_j,s_k)\circ\pi)\frac{\partial}{\partial
y_k}.
\end{equation}
%Hence
%\begin{equation}\label{Lie-1}
%[X^h,Y^h]=[X,Y]^h+\sum_{j,k=1}^3y_j(g(R^D(X,Y)s_j,s_k)\circ\pi)\frac{\partial}{\partial
%y_k}
%\end{equation}
%for every vector fields $X,Y$ on ${\mathscr N}$. Let
%$\sigma\in{\mathcal Z}$. Using the standard identification
%$T_{\sigma}(\Lambda^2_{+}T_{\pi(\sigma)}M)\cong
%\Lambda^2_{+}T_{\pi(\sigma)}M$,  we obtain from (\ref{Lie-1}) the
%well-known formula
%\begin{equation}\label{Lie-2h}
%[X^h,Y^h]_{\sigma}=[X,Y]^h_{\sigma}+R^D_{p}(X,Y)\sigma, \quad
%p=\pi(\sigma).
%\end{equation}
Let $\sigma\in{\mathcal Z}$. Under the standard identification
$T_{\sigma}(\Lambda^2_{+}T_{\pi(\sigma)}M)\cong
\Lambda^2_{+}T_{\pi(\sigma)}M$, formula (\ref{hl}) implies the
well-known identity
\begin{equation}\label{Lie-2h}
[X^h,Y^h]_{\sigma}=[X,Y]^h_{\sigma}+R^D_{p}(X,Y)\sigma, \quad
p=\pi(\sigma).
\end{equation}

For a fixed $\sigma\in{\mathcal Z}$, take an oriented orthonormal
frame $E_1,...,E_4$ such that $(E_3)_p=K_{\sigma}(E_2)_p$,
$(E_4)_p=K_{\sigma}(E_1)_p$, $p=\pi(\sigma)$, so $s_3(p)=\sigma$,
where $s_3$ is defined by means of $E_1,...,E_4$ via
(\ref{s-basis}). Define coordinates $(\tilde x_{\alpha},y_j)$ as
above by means of this frame and a coordinate system of $M$ at $p$.
Set
\begin{equation}\label {V-12}
\begin{array}{c}
V_1=\displaystyle{(1-y_2^2)^{-1/2}(y_3\frac{\partial}{\partial y_1}-y_1\frac{\partial}{\partial y_3})},\\[6pt]
V_2=\displaystyle{(1-y_2^2)^{-1/2}(-y_1y_2\frac{\partial}{\partial
y_1}+(1-y_2^2)\frac{\partial}{\partial
y_2}-y_2y_3\frac{\partial}{\partial y_3})}.
\end{array}
\end{equation}
Then $V_1,V_2$ is a $g$-orthonormal frame of vertical vector fields
in a neighbourhood of the point $\sigma$ such that
$(V_1)_{\sigma}=s_1(p)$, $(V_2)_{\sigma}=s_2(p)$. An easy
computation using (\ref{hl}) gives
\begin{equation}\label{XhVk}
[X^h,V_1]_{\sigma}=g(D_{X}s_1,s_2)s_2(p),\quad
[X^h,V_2]_{\sigma}=-g(D_{X}s_1,s_2)s_1(p).
\end{equation}
Then the Koszul formula for the Levi-Civita connection implies that
the vectors $(\widetilde{D}_{V_k}V_l)_{\sigma}$ are orthogonal to
every horizontal vector $X^h_{\sigma}$, hence they are vertical
vectors. It follows that the fibers of the bundle $\pi:({\mathcal
Z},h_t)\to (M,g)$ are totally geodesic submanifolds.

\smallskip
\noindent {\bf Notation}. We denote the Levi-Civita connection of
$(M,g)$ by $\nabla$.

\smallskip

The Koszul formula, identity (\ref{Lie-2h}) and the fact that the
fibres of the twistor bundle are totally geodesic submanifolds imply
the following formulas.

\begin{lemma}\label{LC} {\rm (\cite{DM})}
If $X,Y$ are vector fields on $M$ and $V$ is a vertical vector field
on ${\mathcal Z}$, then
\begin{equation}\label{D-hh}
(\widetilde
D_{X^h}Y^h)_{\sigma}=(\nabla_{X}Y)^h_{\sigma}+\frac{1}{2}R^D_{\pi(\sigma)}(X,Y)\sigma,
,\quad \sigma\in{\mathcal Z}.
\end{equation}
Also, $\widetilde D_{V}X^h={\mathcal H}\widetilde D_{X^h}V$, where
${\mathcal H}$ means "the horizontal component". Moreover,
\begin{equation}\label{D-vh}
h_t(\widetilde
D_{V}X^h,Y^h)_{\sigma}=-\frac{t}{2}g(R^D_{\pi(\sigma)}(X,Y)\sigma,V),\quad
\sigma\in{\mathcal Z}.
\end{equation}
\end{lemma}

\section{The second fundamental form of an  almost Hermitian structure as a map into the twistor space}

Let $(M,g,J)$ be an almost Hermitian manifold of dimension four.
Define a section ${\frak J}$ of $\Lambda^2TM$ by
$$
g({\frak J},X\wedge Y)=\frac{1}{2}g(JX,Y),\quad X,Y\in TM.
$$
Note that the section $2{\frak J}$ is dual to the fundamental
$2$-form of $(M,g,J)$. Consider $M$ with the orientation induced by
the almost complex structure $J$. Then ${\frak J}$ takes its values
in the (positive) twistor space ${\mathcal Z}$ of the Riemannian
manifold $(M,g)$.

\smallskip
\noindent {\bf Notation}. Denote by $\widetilde{D}^{\ast}$ the
connection on the pull-back bundle ${\frak J}^{\ast}T{\cal Z}$
induced by the Levi-Civita connection $\widetilde{D}$ of $({\mathcal
Z},h_t)$. The connections $D$ on $TM$ and $\widetilde{D}^{\ast}$ on
${\frak J}^{\ast}T{\cal Z}$ induce a connection $\widehat D$ on the
bundle $Hom(TM,{\frak J}^{\ast}T{\mathcal Z})$.

\begin{prop}\label{covder-dif}
For every $p\in M$ and every $X,Y\in T_pM$,
$$
\begin{array}{c}
(\widehat D _{X}{\frak J}_{\ast})(Y)
=\displaystyle{\frac{1}{2}}[D^{2}_{XY}{\frak J} -g(D^{2}_{XY}{\frak
J},{\frak J}){\frak J}(p)
+ D^{2}_{YX}{\frak J}-g(D^{2}_{YX}{\frak J},{\frak J}){\frak J}(p)\\[6pt]
\hfill -D_{T(X,Y)}{\frak J}-(T(X,Y))^h_{{\frak J}(p)}]\\[6pt]
 +(\widetilde{D}_{D_{X}{\mathfrak J}}Y^h)_{{\frak J}(p)} +(\widetilde{D}_{D_{Y}{\mathfrak J}}X^h)_{{\frak
 J}(p)}.
\end{array}
$$
where  $D^{2}_{XY}{\frak J}=D_X D_Y{\frak J}-D_{D_{X}Y}{\frak J}$ is
the second covariant derivative of the section ${\frak J}$ of
$\Lambda^2_{+}TM$ and $T$ is the torsion of $D$.
\end{prop}

\noindent {\bf Proof}. Extend $X$ and $Y$ to vector fields in a
neighbourhood of the point $p$. Take an oriented orthonormal frame
$E_1,...,E_4$ near $p$ such that $E_3=JE_2$, $E_4=JE_1$, so ${\frak
J}=s_3$. Introduce coordinates $(\tilde x_{\alpha},y_j)$ as above by
means of this frame and a coordinate system of $M$ at $p$. Let
$V_1,V_2$ be the $g$-orthonormal frame of vertical vector fields in
a neighbourhood of the point $\sigma={\frak J}(p)$ defined by
(\ref{V-12}). For this frame,  $V_1\circ{\frak J}=s_1$,
$V_2\circ{\frak J}=s_2$. Note also that $[V_1,V_2]_{\sigma}=0$. This
and the Koszul formula imply  $(\widetilde D_{V_k}V_l)_{\sigma}=0$
since  $\widetilde D_{V_k}V_l$ are vertical vector fields,
$k,l=1,2$. Thus  $\widetilde D_{W}V_l=0$, $l=1,2$, for every
vertical vector $W$ at $\sigma$. Considering ${\frak J}$ as a
section of $\Lambda^2_{+}TM$, we have
$$
{\frak J}_{\ast}\circ Y=Y^h\circ{\frak J}+D_{Y}{\frak
J}=Y^h\circ{\frak J}+\sum_{k=1}^2g(D_{Y}{\frak
J},s_k)(V_k\circ{\frak J}),
$$
hence
$$
\begin{array}{c}
\widetilde{D}^{\ast}_{X}({\frak J}_{\ast}\circ Y)=(\widetilde
D_{{\frak J}_{\ast}X}Y^h)\circ{\frak J}
+\sum\limits_{k=1}^2g(D_{Y}{\frak J},s_k)(\widetilde D_{{\frak J}_{\ast}X}V_k)\circ{\frak J}\\[6pt]
+\sum\limits_{k=1}^2[g(D_X D_{Y}{\frak J},s_k)+g(D_{Y}{\frak
J},D_{X}s_k)](V_k\circ{\frak J}).
\end{array}
$$
This, in view of Lemma~\ref{LC},  implies
$$
\begin{array}{c}
\widetilde{D}^{\ast}_{X_p}({\frak J}_{\ast}\circ Y)=(\nabla_{X}Y)^h_{{\frak J}(p)}+\displaystyle{\frac{1}{2}}R^D(X,Y){{\frak J}(p)}\\[6pt]
+(\widetilde{D}_{D_{X}{\mathfrak J}}Y^h)_{{\frak J}(p)} +(\widetilde{D}_{D_{Y}{\mathfrak J}}X^h)_{{\frak J}(p)}\\[6pt]
+\sum\limits_{k=1}^2g(D_{X_p} D_{Y}{\frak
J},s_k(p))s_k(p)+\sum\limits_{k=1}^2[g(D_{Y_p}{\frak
J},s_k(p))[X^h,V_k]_{{\frak
J}(p)}\\[6pt]
+g(D_{Y_p}{\frak J},D_{X_p}s_k)s_k(p)].
\end{array}
$$
By the convention for the curvature tensor,
$$
R^D(X,Y){\frak J}(p)=D_{D_{X_p}Y}{\frak J}-D_{D_{Y_p}X}{\frak
J}-D_{T(X,Y)_p}{\frak J}-D_{X_p}D_{Y}{\frak J}+D_{Y_p}D_{X}{\frak
J}.
$$

 Identities (\ref{XhVk}) imply
$$
\sum_{k=1}^2[g(D_{Y_p}{\frak J},s_k(p))[X^h,V_k]_{{\frak
J}(p)}+g(D_{Y_p}{\frak J},D_{X_p}s_k(p))s_k(p)]=0.
$$
Since $g(D_{X}{\frak J},{\frak J})=0$ for every $X$, we have
$$
g(D_{Y} D_{X}{\frak J},{\frak J})=-g(D_{X}{\frak J},D_{Y}{\frak
J})=g(D_{X} D_{Y}{\frak J},{\frak J}).
$$
Hence
$$
\sum_{k=1}^2g(D_{X_p} D_{Y}{\frak
J},s_k(p))s_k(p)=D_{X_p}D_{Y}{\frak J}-\frac{1}{2}g(D_{X_p}
D_{Y}{\frak J} +D_{Y_p}D_{X}{\frak J},{\frak J}){\frak J}(p).
$$
The identities above give
$$
\begin{array}{c}
\widetilde{D}^{\ast}_{X_p}({\frak J}_{\ast}\circ Y)=(\nabla_{X}Y)^h_{{\frak J}(p)}+\\[6pt]
\displaystyle{\frac{1}{2}[(D_{X_p}D_{Y}{\frak J}+D_{Y_p}D_{X}{\frak
J})-g(D_{X_p} D_{Y}{\frak J} +D_{Y_p}D_{X}{\frak J}){\frak J}(p)}]\\[6pt]
+\displaystyle{\frac{1}{2}[D_{D_{X_p}Y}{\frak J}-D_{D_{Y_p}X}{\frak
J}-D_{T(X,Y)_p}{\frak J}]}%\\[8pt]
+(\widetilde{D}_{D_{X}{\mathfrak J}}Y^h)_{{\frak J}(p)} +(\widetilde{D}_{D_{Y}{\mathfrak J}}X^h)_{{\frak J}(p)}.\\[6pt]
\end{array}
$$
This and (\ref{D-nabla}) imply the desired formula for
$$
(\widehat D_{X}{\frak
J}_{\ast})(Y)=\widetilde{D}^{\ast}_{X_p}({\frak J}_{\ast}\circ Y)
-(D_{X_p}Y)^h_{{\frak J}(p)}-D_{D_{X_p}Y}{\frak J}.
$$

%\hfill $\Box$

Proposition~\ref{covder-dif} and Lemma~\ref{LC} imply
\begin{cor}\label{V-H-tr}
For every $p\in M$,
$$
\begin{array}{c}
{\mathcal V}(Trace_g \widehat D{\frak
J}_{\ast})_{p}=Trace_g\{T_{p}M\ni X\to {\mathcal V} D^2_{XX}{\frak
J}\}\\[6pt]
{\mathcal H}(Trace_g \widehat D{\frak
J}_{\ast})_{p}=2Trace_g\{T_{p}M\ni X\to
(\widetilde{D}_{D_{X}{\mathfrak J}}X^h)_{{\frak J}(p)}\},
\end{array}
$$
where ${\mathcal V}$ and ${\mathcal H}$ mean the horizontal and the
vertical component, respectively.
\end{cor}

Note that
$$
g(D^2_{XY}{\frak J},Z\wedge U)=\frac{1}{2}g((D^2_{XY}J)(Z),U),\quad
X,Y,Z,U\in TM.
$$

Every vector of the vertical space ${\mathcal V}_{{\frak J}(p)}$,
$p\in M$, is a linear combination of vectors of the type $Z\wedge
U-JZ\wedge JU$, $Z,U\in T_pM$. Moreover, by (\ref{cs}), each vector
$Z\wedge U-JZ\wedge JU$ is orthogonal to every $a\in\Lambda^2_{-}TM$
since the endomorphisms $J$ and $K_a$  of $T_pM$ commute. Also,
$Z\wedge U-JZ\wedge JU$ is orthogonal to ${\frak J}(p)$. Hence
$Z\wedge U-JZ\wedge JU\in {\mathcal V}_{{\frak J}(p)}$.

%\begin{cor}\label{D}
%$$
%\begin{array}{c}
%2g((D^2_{XX}{\mathfrak J},Z\wedge  U-JZ\wedge JU)=g((D^2_{XX}J)(Z),U)-g(D^2_{XX})(JZ),JU)\\[8pt]
%=g((\nabla^2_{XX}J)(Z),U)-g(\nabla^2_{XX}J)(JZ),JU)\\[8pt]
%+(\nabla_{X}{\mathcal T})(X,JZ,U)+(\nabla_{X}{\mathcal
%T})(X,Z,JU)\\[8pt]
%+ {\mathcal T}(X,Z,(\nabla_{X} J)(U)) -{\mathcal T}(X,U,(\nabla_{X}J)(Z))\\[8pt]
%-{\mathcal T}(X,JZ,(\nabla_{X} J)(JU))+{\mathcal
%T}(X,JU,(\nabla_{X}J)(JZ)).
%\end{array}
%$$
%\end{cor}

\smallskip

\noindent {\bf Assumption}.  Suppose that the metric connection $D$
has skew-symmetric torsion $T$, so ${\mathcal T}(X,Y,Z)=g(T(X,Y),Z)$
is a skew-symmetric $3$-form.

\smallskip

\noindent {\bf Notation}. It is convenient to set
$$
\begin{array}{c}
Q(X,Z\wedge U-JZ\wedge JU)= {\mathcal T}(X,Z,(\nabla_{X} J)(U)) -{\mathcal T}(X,U,(\nabla_{X}J)(Z))\\[6pt]
-{\mathcal T}(X,JZ,(\nabla_{X} J)(JU))+{\mathcal
T}(X,JU,(\nabla_{X}J)(JZ))
\end{array}
$$

\begin{lemma}\label{D}
For every tangent vectors  $X,Y,Z,U\in T_pM$.
$$
\begin{array}{c}
2g((D^2_{XX}{\mathfrak J},Z\wedge  U-JZ\wedge
JU)=2g((\nabla^2_{XX}{\mathfrak J},Z\wedge  U-JZ\wedge JU)\\[6pt]
+(\nabla_{X}{\mathcal T})(X\wedge (Z\wedge JU+JZ\wedge
U))+Q(X,Z\wedge U-JZ\wedge JU)\\[6pt]
-\displaystyle{\frac{1}{2}}{\mathcal
T}(X,JT(X,Z),U)+\displaystyle{\frac{1}{2}}{\mathcal
T}(X,JT(X,JZ),JU).
\end{array}
$$
\end{lemma}

\noindent {\bf Proof}. Extend $X,Y,Z,U$ to vector fields. Then an
easy computation gives
$$
\begin{array}{c}
2g((D^2_{XX}{\mathfrak J},Z\wedge  U-JZ\wedge
JU)=2g((\nabla^2_{XX}{\mathfrak J},Z\wedge  U-JZ\wedge JU)\\[6pt]
-g(\nabla_{X}JT(X,Z),U)-g(T(X,J\nabla_{X}Z),U)-\displaystyle{\frac{1}{2}}g(T(X,JT(X,Z),U)\\[6pt]
+g(\nabla_{X}JT(X,JZ),JU)+g(T(X,J\nabla_{X}JZ),JU)+\displaystyle{\frac{1}{2}}g(T(X,JT(X,JZ),JU)\\[6pt]
-g(T(\nabla_{X}X,JZ),U)-g(T(\nabla_{X}X,Z),JU)\\[6pt]
=-g((\nabla_{X}J)(T(X,Z)),U)+g(\nabla_{X}T)(X,Z),JU)+g(T(X,(\nabla_{X}J)(Z),U)\\[6pt]
+g((\nabla_{X}J)(T(X,JZ)),JU)+g((\nabla_{X}T)(X,JZ),U)-g(T(X,(\nabla_{X}J)(JZ),JU)\\[6pt]
-\displaystyle{\frac{1}{2}}g(T(X,JT(X,Z),U)+\displaystyle{\frac{1}{2}}g(T(X,JT(X,JZ),JU)
\end{array}
$$

This proves the lemma.

%\hfill $\Box$

We also need some properties of the curvature of the connection $D$.

Denote by $R^{\nabla}$ the curvature tensor of the Levi-Civita
connections $\nabla$.

A simple computation gives
\begin{equation}\label{curv-D}
\begin{array}{c}
g(R^D(X,Y)Z,U)=g(R^{\nabla}(X,Y)Z,U)\\[6pt]
-\displaystyle{\frac{1}{2}}\big[(\nabla_{X}{\mathcal
T})(Y,Z,U)-(\nabla_{Y}{\mathcal T})(X,Z,U)\big]\\
+\displaystyle{\frac{1}{4}}\sum\limits_{i=1}^4\big[{\mathcal
T}(X,U,E_i){\mathcal T}(Y,Z,E_i)-{\mathcal T}(X,Z,E_i){\mathcal
T}(Y,U,E_i)\big],
\end{array}
\end{equation}
where $X,Y,Z,U\in T_pM$ and $\{E_1,...,E_4\}$ is an orthonormal
basis of $T_pM$. This implies
\begin{equation}\label{Z-U}
g(R^D(X,Y)Z,U)=-g(R^D(X,Y)U,Z)
\end{equation}

\smallskip

\noindent {\bf Notation}. The $1$-form $\ast{\mathcal T}$ will be
denoted by $\tau$.

\smallskip

Clearly, the form $\tau$ uniquely determines the $3$-form ${\mathcal
T}$, hence the connection $D$.

\smallskip

\noindent {\bf Convention}. We identify $TM$ and $T^{\ast}M$ be
means of the metric $g$ and extend this isomorphism to
identification of $\Lambda^kTM$ and $\Lambda^kT^{\ast}M$. The
exterior product on  $\Lambda^kT^{\ast}M$ is so that
$(\alpha_1\wedge...\wedge\alpha_k)(v_1,...,v_k)=det[\alpha_i(v_j)]$
for $\alpha_i\in T^{\ast}M$ and $v_j\in TM$.

\smallskip

 For a given orthonormal frame $E_1,...,E_4$, it is
convenient to set
$$
E_{ijk}=E_i\wedge E_j\wedge E_k, \quad {\mathcal T}_{ijk}={\mathcal
T}(E_{ijk}).
$$
We consider $\Lambda^3TM$ with the metric for which $E_{ijk}$,
$1\leq i<j<k\leq 4$, is an orthonormal basis.  Then we have
\begin{equation}\label{T-tau}
\begin{array}{c}
{\mathcal T}={\mathcal T}_{123}E_{123}+{\mathcal
T}_{124}E_{124}+{\mathcal T}_{134}E_{134}+{\mathcal
T}_{234}E_{234}\\[6pt]
\tau=-{\mathcal T}_{234}E_1+{\mathcal T}_{134}E_2-{\mathcal
T}_{124}E_3+{\mathcal T}_{123}E_4.
\end{array}
\end{equation}
It follows that
\begin{equation}\label{nabla-T-tau}
\begin{array}{c}
(\nabla_{X}{\mathcal T})(E_{123})=(\nabla_{X}\tau)(E_4),\quad
(\nabla_{X}{\mathcal T})(E_{124})=-(\nabla_{X}\tau)(E_3),\\[6pt]
(\nabla_{X}{\mathcal T})(E_{134})=(\nabla_{X}\tau)(E_2),\quad
(\nabla_{X}{\mathcal T})(E_{234})=-(\nabla_{X}\tau)(E_1).
\end{array}
\end{equation}
The latter identities imply that for every $a\in\Lambda^2_{+}TM$
\begin{equation}\label{d-delta}
d\tau(a)=-\delta{\mathcal T}(a).
\end{equation}

\noindent {\bf Notation}. The Ricci tensor and $\ast$-Ricci tensor
are defined by
$$
\begin{array}{c}
\rho_D(X,Y)=Trace\{Z\to g(R^D(X,Z)Y,Z)\},\\[6pt]
\rho_{D}^{\ast}(X,Y)=Trace\{Z\to g(R^D(JZ,X)JY,Z)\}.
\end{array}
$$

Note that $\rho^{\ast}_{\nabla}(X,Y)=\rho^{\ast}_{\nabla}(JY,JX)$ by
the identity
$g(R^{\nabla}(JZ,X)JY,Z)=\\-g(R^{\nabla}(Z,JY)J(JX),Z)$.
Equivalently, $\rho^{\ast}_{\nabla}(X,JX)=0$.
\smallskip
\begin{prop}\label {r,star-r}
$$
\begin{array}{c}

\rho_D(X,Y)=\rho_{\nabla}(X,Y)-\displaystyle{\frac{1}{2}}\delta{\mathcal
T}(X,Y)-2g(\imath_{X}{\mathcal T},\imath_{Y}{\mathcal
T})\\[6pt]

\rho^{\ast}_D(X,Y)=\rho_{\nabla}^{\ast}(X,Y)+d{\mathcal T}({\frak
J}\wedge X\wedge JY)+(\nabla_{JY}{\mathcal T})({\frak
J}\wedge X)\\[6pt]
-\displaystyle{\frac{1}{4}}[2g(T({\frak J}),T(X\wedge
JY))+\chi(X,Y)],
\end{array}
$$
where $\chi(X,Y)=Trace \{\Lambda^2_{+}T_pM\ni a\to {\mathcal
T}(a\wedge X){\mathcal T}(({\frak J}\times a)\wedge JY)\}$. Moreover
$$
\chi(X,Y)=\chi(Y,X),\quad \chi(X,Y)=\chi(JY,JX).
$$
\end{prop}

\noindent {\bf Proof}. The first formula follows directly from
(\ref{curv-D}).

Let $E_1,...,E_4$ be an oriented orthonormal basis of a tangent
space $T_pM$ such that $E_2=JE_1$, $E_4=JE_3$. By (\ref{curv-D})
$$
\begin{array}{c}
\rho^{\ast}_D(X,Y)=\rho_{\nabla}^{\ast}(X,Y)-\displaystyle{\frac{1}{2}}\sum\limits_{k=1}^4[(\nabla_{JE_k}{\mathcal
T})(X,JY,E_k)-(\nabla_{X}{\mathcal T})(JE_k,JY,E_k)]\\
+\displaystyle{\frac{1}{4}}\sum\limits_{i,k=1}^4[{\mathcal
T}(JE_k,E_k,E_i){\mathcal T}(X,JY,E_i)-{\mathcal
T}(JE_k,JY,E_i){\mathcal T}(X,E_k,E_i)].
\end{array}
$$
To compute the second summand in the right-hand side, we apply the
identity
$$
\begin{array}{c}
d{\mathcal T}(JE_k,X,JY,E_k)=(\nabla_{JE_k}{\mathcal
T})(X,JY,E_k)-(\nabla_{X}{\mathcal T})(JE_k,JY,E_k)\\[6pt]
\hfill +(\nabla_{JY}{\mathcal T})(JE_k,X,E_k)-(\nabla_{E_k}{\mathcal
T}(JE_k,X,JY).
\end{array}
$$
Summing up, we get
$$
\begin{array}{c}
-2d{\mathcal T}({\frak J}\wedge X\wedge
JY)=\sum\limits_{k=1}^4[(\nabla_{JE_k}{\mathcal
T})(X,JY,E_k)-(\nabla_{X}{\mathcal T})(JE_k,JY,E_k)]\\[6pt]
+2(\nabla_{JY}{\mathcal T})({\frak J}\wedge X).
\end{array}
$$
Next,
$$
\sum\limits_{i,k=1}^4{\mathcal T}(JE_k,E_k,E_i){\mathcal
T}(X,JY,E_i)=-2g(T({\frak J}),T(X\wedge JY)).
$$
Finally, a direct computation gives
$$
\begin{array}{c}
\sum\limits_{i,k=1}^4{\mathcal T}(JE_k,JY,E_i){\mathcal
T}(X,E_k,E_i)\\[6pt]
=-{\mathcal T}(s_2\wedge X){\mathcal T}(s_3\wedge
JY)+{\mathcal T}(s_3\wedge X){\mathcal T}(s_2\wedge JY)\\[6pt]
=-Trace \{\Lambda^2_{+}T_pM\ni a\to {\mathcal T}(a\wedge X){\mathcal
T}(({\frak J}\times a)\wedge JY)\}=-\chi(X,Y).
\end{array}
$$
We have $\chi(X,Y)=\chi(JY,JX)$ since
$$\sum\limits_{k=1}^n{\mathcal T}(JE_k,JY,E_i){\mathcal T}(X,E_k,E_i)=-\sum\limits_{k=1}^n{\mathcal
T}(E_k,J(JX),E_i){\mathcal T}(JY,JE_k,E_i).
$$
A direct computation shows that $\chi(E_i,E_j)=\chi(E_j,E_i)$,
$i,j=1,...,4$.

\section{Harmonicity of ${\frak J}$ in the case of a Hermitian structure}

Let $\Omega(X,Y)=g(JX,Y)$ be the fundamental $2$-form of the almost
Hermitian manifold $(M,g,J)$. Denote by $N$ the Nijenhuis tensor of
$J$:
$$
N(Y,Z)=-[Y,Z]+[JY,JZ]-J[Y,JZ]-J[JY,Z].
$$
Clearly $ N(Y,Z)=-N(Z,Y), \quad N(JX,Y)=N(X,JY)=-JN(X,Y)$.

It is well-known (and easy to check) that
\begin{equation}\label{nJ}
2g((\nabla_XJ)(Y),Z)=d\Omega(X,Y,Z)-d\Omega(X,JY,JZ)+g(N(Y,Z),JX),
\end{equation}
for all $X,Y,Z\in TM$.  Note also that
$$
g(\nabla_{X}J)(Y),Z)=(\nabla_{X}\Omega)(Y,Z).
$$

{\it Suppose that the almost complex structure $J$ is integrable}.
This is equivalent to
\begin{equation}\label{nablaJXJY}
(\nabla_{X}J)(Y)=(\nabla_{JX}J)(JY),\quad X,Y\in TM,
\end{equation}
\cite[Corollary 4.2]{G}. Let $B$ be the vector field on $M$ dual to
the Lee form $\theta=-\delta\Omega\circ J$ with respect to the
metric $g$. Then (\ref{nJ}) and the identity
$d\Omega=\theta\wedge\Omega$ imply the following well-known formula
\begin{equation}\label{nablaJ}
2(\nabla_XJ)(Y)=g(JX,Y)B-g(B,Y)JX+g(X,Y)JB-g(JB,Y)X.
\end{equation}
Since $\nabla$ and $D$ are metric connections,
$$
g(\nabla_X{\frak J},Y\wedge
Z)=\displaystyle{\frac{1}{2}}g((\nabla_XJ)(Y),Z),\quad g(D_X{\frak
J},Y\wedge Z)=\displaystyle{\frac{1}{2}}g((D_XJ)(Y),Z)
$$
Hence
\begin{equation}\label{nabla-frak}
\nabla_X{\frak J}=\displaystyle{\frac{1}{2}}(JX\wedge B+X\wedge JB).
\end{equation}
In order  to compute the vertical part of $(Trace_g \widehat D{\frak
J}_{\ast})_{p}$, we apply Corollary~\ref{V-H-tr}  and Lemma~\ref{D}.
First, identities (\ref{nablaJ}) and (\ref{nabla-frak})  imply
(\cite{DHM})
$$
\begin{array}{r}
(\nabla^2_{XY}J)(Z)=\displaystyle{\frac{1}{2}}\big[g((\nabla_{X}J)(Y),Z)B-g(B,Z)(\nabla_{X}J)(Y)\\[6pt]
g(Y,Z)(\nabla_{Y}J)(B)-g((\nabla_{Y}J)(B),Z)Y\\[6pt]

+g(JY,Z)\nabla_{X}B-g(\nabla_{X}B,Z)JY+g(Y,Z)J\nabla_{X}B-g(J\nabla_{X}B,Z)Y\big].
\end{array}
$$

It follows that
\begin{equation}\label{nablaXX}
2g(Trace\{X\to(\nabla^2_{XX}J)(Z)\},U)=||B||^2g(Z,JU)-d\theta(JZ,U)-d\theta(Z,JU).
\end{equation}

Let $E_1,...,E_4$ be an orthonormal basis of a tangent space $T_pM$
with  $E_2=JE_1$, $E_4=JE_3$.  Define $s_1,s_2,s_3$ by means of this
basis via (\ref{s-basis}). Thus, ${\frak J}=s_1$ and $s_2,s_3$ is a
$g$-orthonormal basis of the vertical space ${\mathcal V}_{{\frak
J}(p)}$.

Identity (\ref{nablaXX}) implies
$$
\begin{array}{c}
g(Trace\{(\nabla^2_{XX}J)(E_1)\},E_3) +
g(Trace\{(\nabla^2_{XX}J)(E_4)\},E_2) = -d\theta(s_3), \\[6pt]
g(Trace\{(\nabla^2_{XX}J)(E_1)\},E_4) +
g(Trace\{(\nabla^2_{XX}J)(E_2)\},E_3) =d\theta(s_2)
\end{array}
$$

Also, it follows from (\ref{d-delta})
$$
Trace\{(\nabla_{X}{\mathcal T})(X,s_2)\}= d\tau(s_2),\quad
Trace\{(\nabla_{X}{\mathcal T})(X,s_3)\}=d\tau(s_3).
$$
Applying (\ref{nablaJ}), we obtain
$$
\begin{array}{c}
Trace \{Q(X,s_2)\}=\sum\limits_{k=1}^4Q(E_k,E_1\wedge E_3 -E_2\wedge E_4)\\[6pt]
=-g(B,E_1){\mathcal T}_{123}+g(B,E_2){\mathcal
T}_{124}+g(B,E_3){\mathcal T}_{134}-g(B,E_4){\mathcal
T}_{234}\\[6pt]
-{\mathcal T}(E_1,E_3,JB)+{\mathcal T}(E_1,E_4,B)+{\mathcal T}(E_2,E_3,B)+{\mathcal T}(E_2,E_4,JB)\\[6pt]
=g(B,E_1){\mathcal T}_{123}-g(B,E_3){\mathcal
T}_{134}+g(B,E_2){\mathcal T}_{234}-g(B,E_4){\mathcal
T}_{124}\\[6pt]
=g(T(E_1\wedge E_4+E_2\wedge E_3),B)=g(T(s_3),B).
\end{array}
$$
Also,
$$
\begin{array}{c}
Trace\{X\to -{\mathcal T}(X,JT(X,E_1,),E_3)+{\mathcal
T}(X,JT(X,E_2),E_4)\}\\[6pt]
=-T_{213}T_{243}+T_{412}T_{413}-T_{124}T_{134}+T_{321}T_{324}=0.
\end{array}
$$
Similarly,
$$
\begin{array}{c}
Trace\{Q(X,s_3)\}=\sum\limits_{k=1}^4Q(E_k,E_1\wedge E_4 +E_2\wedge E_3)\\[6pt]
=-g(T(E_1\wedge E_3+E_4\wedge E_2),B)=-g(T(s_2),B)
\end{array}
$$
and
$$
Trace\{X\to -{\mathcal T}(X,JT(X,E_1,),E_4)-{\mathcal
T}(X,JT(X,E_2),E_3)\}=0.
$$
Thus, by Lemma~\ref{D},
$$
\begin{array}{c} 2g(Trace\,D^2_{XX}{\frak
J},s_2)=-d\theta(s_3)+d\tau(s_3)+g(T(s_3),B)\\[6pt]
2g(Trace\,D^2_{XX}{\frak
J},s_3)=d\theta(s_2)-d\tau(s_2)-g(T(s_2),B).
\end{array}
$$
Finally, note that for every $X,Y\in TM$, the $2$-vector $X\wedge
Y-JX\wedge JY$ is a linear combination of $s_2$ and $s_3$. Thus, by
Corollary~\ref{V-H-tr}, ${\mathcal V}Trace_g \widehat D{\frak
J}_{\ast}=0$ if and only if the $2$-form
$d\theta-d\tau-\imath_{B}{\mathcal T}$ is of type $(1,1)$ with
respect to $J$.

 Now, suppose  $B_p\neq 0$ at a point $p\in M$. Take an
oriented orthonormal basis $E_1,...,E_4$ of $T_pM$ such that
$E_2=JE_1$, $E_3=\frac{B_p}{||B_p||}$, $E_4=JE_3$. Using this basis,
define $s_1,s_2,s_3$ via (\ref{s-basis}) so that ${\frak J}=s_1$. By
(\ref{D-vh}),
$$
h_t(Trace\{T_pM\ni X\to D_{D_{X}{\frak J}}X^h\},Z^h)_{{\frak
J}(p)}=-\frac{t}{2}\sum\limits_{k=1}^4g(R^D(E_k,Z){\frak
J}(p),D_{E_k}{\frak J}).
$$
for every $Z\in T_pM$.  Note that
$$
(D_{X}J)(Y)=(\nabla_{X}J)(Y)+\frac{1}{2}[T(X,JY)-JT(X,Y)].
$$
This identity implies
$$
\begin{array}{c}
D_{X}{\frak J}=g(D_{X}{\frak J},s_2)s_2+g(D_{X}{\frak
J},s_3)s_3\\[6pt]
=\displaystyle{\frac{1}{2}}[g((D_{X}J)(E_1),E_3)+g((D_{X}J)(E_4),E_2)]s_2\\[6pt]
+\displaystyle{\frac{1}{2}}[g((D_{X}J)(E_1),E_4)+g((D_{X}J)(E_2),E_3)]s_3\\[6pt]
=\nabla_{X}{\frak J}+\displaystyle{\frac{1}{2}}[{\mathcal
T}(X,s_3)s_2-{\mathcal T}(X,s_2)s_3]
\end{array}
$$
Thus, by (\ref{nabla-frak}),
$$
D_{X}{\frak J}=\displaystyle{\frac{1}{2}}[JX\wedge B+X\wedge
JB+{\mathcal T}(X,s_3)s_2-{\mathcal T}(X,s_2)s_3].
$$
Then, applying Lemma~\ref{CurOp}, we compute:
$$
\begin{array}{c}
\sum\limits_{k=1}^4g(R^D(E_k,Z){\frak J}(p),D_{E_k}{\frak J})\\[6pt]
=-\displaystyle{\frac{1}{2}}||B_p||[g({\mathcal R}^D(E_1\wedge
Z),s_2)+g({\mathcal
R}^D(E_2\wedge Z),s_3)]\\[6pt]
+\displaystyle{\frac{1}{2}}\sum\limits_{k=1}^4[{\mathcal
T}(E_k,s_2)g({\mathcal R}^D(E_k\wedge Z),s_2)+{\mathcal
T}(E_k,s_3)g({\mathcal R}^D(E_k\wedge Z),s_3)]\\[6pt]
=-\displaystyle{\frac{1}{2}}\big[||B_p|||g({\mathcal R}^D(E_1\wedge
Z),s_2)+||B_p||g({\mathcal
R}^D(E_2\wedge Z),s_3)\\[6pt]
\hfill - g({\mathcal R}^D(T(s_2)\wedge Z),s_2)-g({\mathcal
R}^D(T(s_3)\wedge Z),s_3)\big].
\end{array}
$$
Using (\ref{Z-U}), it is easy to see that
$$
||B_p||[g({\mathcal R}^D(E_1\wedge Z),s_2)+g({\mathcal
R}^D(E_2\wedge Z),s_3)]=\rho_D(Z,B)-\rho^{\ast}_D(Z,B).
$$
Next,
$$
\begin{array}{c}
T(s_2)=-{\mathcal T}_{124}E_1-{\mathcal T}_{123}E_2+{\mathcal
T}_{234}E_3+{\mathcal T}_{134}E_4,\\[6pt]
T(s_3)={\mathcal T}_{123}E_1-{\mathcal T}_{124}E_2-{\mathcal
T}_{134}E_3+{\mathcal T}_{234}E_4.
\end{array}
$$
Then a direct computation gives
\begin{equation}\label{T(s)}
g({\mathcal R}^D(T(s_2)\wedge Z),s_2)+g({\mathcal R}^D(T(s_3)\wedge
Z),s_3)=\rho_{D}(Z,\tau)-\rho_{D}^{\ast}(Z,\tau).
\end{equation}
Thus, if $B_p\neq 0$,
\begin{equation}\label{Tr-h}
\begin{array}{c}
h_t(Trace\{T_pM\ni X\to D_{D_{X}{\frak J}}X^h\},Z^h)_{{\frak
J}(p)}\\
=\displaystyle{\frac{t}{4}}[\rho_D(Z,B)-\rho^{\ast}_D(Z,B)-\rho_{D}(Z,\tau)+\rho_{D}^{\ast}(Z,\tau)].
\end{array}
\end{equation}
If $B_p=0$ for a point $p\in M$, then $\nabla_{X}{\frak J}=0$ for
every $X\in T_pM$ by (\ref{nabla-frak}). Now, the computation above
holds for every oriented orthonormal basis $E_1,...,E_4$ of $T_pM$
with $E_2=JE_1$ and $E_4=JE_3$, and we get identity (\ref{Tr-h})
with vanishing first two summands on right hand side. Hence
${\mathcal H}Trace_g \widehat D{\frak J}_{\ast}=0$ if and only if
for every tangent vector $Z\in TM$
$\rho_D(Z,B)-\rho^{\ast}_D(Z,B)-\rho_{D}(Z,\tau)+\rho_{D}^{\ast}(Z,\tau)=0$.

\smallskip

By Corollary~\ref{har-D}, the considerations above prove the
following.
\begin{theorem}\label{har-Int}
An integrable almost Hermitian structure $J$ determines a harmonic
map from $(M,g)$ into $({\mathcal Z},h_t)$ if and only if $2$-form
$d\theta-d\tau-\imath_{B}{\mathcal T}$ is of type $(1,1)$ with
respect to $J$  and for every $Z\in TM$
$$
\rho_D(Z,B)-\rho^{\ast}_D(Z,B)-\rho_{D}(Z,\tau)+\rho_{D}^{\ast}(Z,\tau)=0.
$$
\end{theorem}

Let $D$ be the metric connection with totally skew-symmetric torsion
determined by the $1$-form $\tau=\theta$. Then the conditions in
Theorem~\ref{har-Int} are trivially satisfied.

In order to specify the torsion $3$-form ${\mathcal T}=-\star
\theta$, recall that
$$
\theta=(\ast\,d\,\ast\Omega)\circ J=(\ast\,d\Omega)\circ J.
$$
Hence $\ast(\theta\circ J)=-\ast^2 d\Omega=d\Omega$. For any
$1$-form $\alpha$,
$$
(\ast\alpha)(X,Y,Z)=-(\ast(\alpha\circ J))(JX,JY,JZ).
$$
Then
$
{\mathcal T}(X,Y,Z)=(\ast(\alpha\circ
J))(JX,JY,JZ)=d\Omega(JX,JY,JZ).
$
The latter identity, (\ref{nablaJXJY}) and (\ref{nJ}) imply
$(D_{X}J)(Y)=0$. Thus $D$ is the Bismut-Strominger connection. Of
course, the fact that ${\mathfrak J}$ is harmonic when we endow $M$
with the Bismut-Strominger connection follows directly from
Corollary~\ref{V-H-tr}.

\smallskip
 \noindent {\bf Example 1}. It has been observed in
\cite{Tr} that every Inoue surface $M$ of type $S^0$ admits a
locally conformal K\"ahler metric $g$  for which the Lee form
$\theta$ is nowhere vanishing, see below. %(see also \cite{DO}).
Define the metric $h_t$ by means of the Levi-Civita connection of
the metric $g$. It is shown in \cite{DHM} that in this case the map
${\mathfrak J}: (M,g)\to ({\mathcal Z},h_t)$ is not harmonic. But,
${\mathfrak J}$ is harmonic if we use the Bismut-Strominger
connection to define $h_t$.  It is natural to ask if there are other
metric connections with skew torsion on $M$ for which ${\mathfrak
J}$ is harmonic.

First,  recall  the construction of the Inoue surfaces of type $S^0$
(\cite{Inoue}). Let $A\in SL(3,\mathbb{Z})$ be a matrix with a real
eigenvalue $\alpha > 1$ and two complex eigenvalues $\beta$ and
$\overline{\beta}$, $\beta\neq\overline{\beta}$. Choose eigenvectors
$(a_1,a_2,a_3)\in{\mathbb R}^3$ and $(b_1,b_2,b_3)\in {\mathbb C}^3$
of $A$ corresponding to $\alpha$ and $\beta$, respectively. Then the
vectors $(a_1,a_2,a_3), (b_1,b_2,b_3),
(\overline{b_1},\overline{b_2},\overline{b_3})$  are
$\mathbb{C}$-linearly independent. Denote the upper-half plane in
$\mathbb{C}$ by ${\bf H}$ and let $\Gamma$ be the group of
holomorphic automorphisms of ${\bf H}\times {\mathbb C}$ generated
by
$$g_o:(w,z)\to (\alpha w,\beta z), \quad g_i:(w,z)\to (w+a_i,z+b_i), \>i=1,2,3 .$$
The group $\Gamma$ acts on ${\bf H}\times{\mathbb C}$ freely and
properly discontinuously.  Then $M=({\bf H}\times {\mathbb
C})/\Gamma$ is a complex surface known as Inoue surface of type
$S^0$.

Following \cite{Tr}, consider on ${\bf H}\times {\mathbb C}$ the
Hermitian metric
$$
g=\frac{1}{v^2}(du\otimes du+dv\otimes dv)+v(dx\otimes dx+dy\otimes
dy),\quad u+iv\in{\bf H}, \quad x+iy\in {\mathbb C}.
$$
This metric is invariant under the action of the group $\Gamma$, so
it descends to a Hermitian metric on $M$ which we denote again by
$g$. Instead on $M$, we  work with $\Gamma$-invariant objects on
${\bf H}\times{\mathbb C}$.  Let $\Omega$ be the fundamental
$2$-form of the Hermitian structure $(g,J)$ on ${\bf
H}\times{\mathbb C}$, $J$ being the standard complex structure. Then
$$ d\Omega=\frac{1}{v}dv\wedge\Omega. $$ Hence the Lee form is
$\theta=d\ln v$. In particular, $d\theta=0$, i.e. $(g,J)$ is a
locally conformal K\"ahler structure. Set
$$
E_1=v\frac{\partial}{\partial u},\quad E_2=v\frac{\partial}{\partial
v},\quad E_3=\frac{1}{\sqrt v}\frac{\partial}{\partial x},\quad
E_4=\frac{1}{\sqrt v}\frac{\partial}{\partial y}.
$$
These are $\Gamma$-invariant vector fields constituting an
orthonormal basis such that $JE_1=E_2$, $JE_3=E_4$. Note that the
vector field dual to the Lee form is $B=E_2$. The non-zero Lie
brackets of $E_1,...,E_4$ are
$$
[E_1,E_2]=-E_1,\quad [E_2,E_3]=-\frac{1}{2}E_3,\quad
[E_2,E_4]=-\frac{1}{2}E_4.
$$
Then we have the following table for the Levi-Civita connection
$\nabla$ of $g$ (\cite{DHM}):
\begin{equation}\label{LC-S0}
\begin{array}{c}
\nabla_{E_1}E_1=E_2,\quad \nabla_{E_1}E_2=-E_1,\\[6pt]
\displaystyle{\nabla_{E_3}E_2=\frac{1}{2}E_3,\quad
\nabla_{E_3}E_3=-\frac{1}{2}E_2, \quad
\nabla_{E_4}E_2=\frac{1}{2}E_4,\quad
\nabla_{E_4}E_4=-\frac{1}{2}E_2}\\[6pt]
{\rm{all~ other}}~ \nabla_{E_i}E_j=0.
\end{array}
\end{equation}

Denote  the dual basis of $E_1,...,E_4$ by $\eta_1,...,\eta_4$. Any
$\Gamma$-invariant $1$-form $\tau$ is of the form $
\tau=a_1\eta_1+...+a_4\eta_4$ where $a_1,...,a_4$ are real
constants. Denote by $D$ the connection which is metric w.r.t. $g$
and has skew-symmetric torsion determined by $\tau$. Note that the
$3$-form ${\mathcal T}$ corresponding to $\tau$ is $ {\mathcal
T}=a_4E_{123}-a_3E_{124}+a_2E_{134}-a_1E_{234}. $

We have   $d\eta_1=\eta_1\wedge\eta_2$, $d\eta_2=0$,
$d\eta_3=\frac{1}{2}\eta_2\wedge \eta_3$,
$d\eta_4=\frac{1}{2}\eta_2\wedge \eta_4$. It follows that
$d\theta-d\tau-\imath_{B}{\mathcal T}$ is of type $(1,1)$ exactly
when $a_3=a_4=0$. Hence $
\rho_D(Z,B)-\rho^{\ast}_D(Z,B)-\rho_{D}(Z,\tau)+\rho_{D}^{\ast}(Z,\tau)=0
$ for every $Z$ exactly when
\begin{equation}\label{cur-eq}
-a_1\rho_D(Z,E_1)+(1-a_2)\rho_D(Z,E_2)+a_1\rho^{\ast}_D(Z,E_1)-(1-a_2)\rho^{\ast}_D(Z,E_2)=0,
\end{equation}
The torsion $T$  corresponding to $\tau$ is given by
$$
\begin{array}{c}
T(E_1,E_2)=0\quad T(E_1,E_3)=a_2E_4,\quad
T(E_1,E_4)=-a_2E_3,\\[6pt]

T(E_2,E_3)=-a_1E_4,\quad T(E_2,E_4)=a_1E_3,\quad
T(E_3,E_4)=a_2E_1-a_1E_2.
\end{array}
$$
Then, using (\ref{curv-D}) and (\ref{LC-S0}), we get the following
table after a tedious computation.
$$
\begin{array}{c}

\rho_D(E_1,E_1)=-\frac{1}{2}a_2^2,\quad
\rho_D(E_2,E_1)=\frac{1}{2}a_1a_2,\quad
\rho_D(E_3,E_1)=\rho_D(E_4,E_1)=0 \\[6pt]

 \rho_D(E_1,E_2)=\frac{1}{2}a_1a_2,\quad
\rho_D(E_2,E_2)=-\frac{3}{2}-\frac{1}{2}a_1^2,\quad
\rho_D(E_3,E_2)=\rho_D(E_4,E_2)=0\\[8pt]

~\rho^{\ast}_D(E_1,E_1)=-1+\frac{1}{2}a_2,\quad
\rho^{\ast}_D(E_2,E_1)=\rho^{\ast}_D(E_3,E_1)=\rho^{\ast}_D(E_4,E_1)=0\\[6pt]

\rho^{\ast}_D(E_1,E_2)=-a_1,\quad
~\rho^{\ast}_D(E_2,E_2)=-1-\frac{1}{2}a_2,\quad
\rho^{\ast}_D(E_3,E_2)=\rho^{\ast}_D(E_4,E_2)=0.

\end{array}
$$
It follows that identity (\ref{cur-eq}) holds for every $Z$ only for
$a_1=0$, $a_2=1$, i.e. $\tau=\theta$. Thus, the complex structure of
an Inoue surface $M$ of type $S^0$  is a harmonic map from $(M,g)$
into its twistor space $({\mathcal Z}, h_t^D)$  only when $D$ is the
Bismut-Strominger connection.

\smallskip
 \noindent {\bf Example 2}. Recall that a primary Kodaira surface $M$ is the quotient of
${\mathbb C}^2$ by  a group of transformations acting freely and
properly discontinuously \cite[p.787]{Kodaira}. This group is
generated by the affine transformations
$\varphi_k(z,w)=(z+a_k,w+\overline{a}_kz+b_k)$ where $a_k$, $b _k$,
$k=1,2,3,4$, are complex numbers such that $ a_1=a_2=0$, $b_2\neq
0$, $Im(a_3{\overline a}_4)=m b_1\neq 0$ for some integer $m>0$. The
quotient space is compact.

It is well-known that $M$ can also be described  as the quotient of
${\mathbb C}^2$ endowed with a group structure by a discrete
subgroup $\Gamma$. The multiplication on ${\Bbb C}^2$ is defined by
$$
(a,b).(z,w)=(z+a,w+\overline{a}z+b),\quad (a,b), (z,w)\in  {\Bbb
C}^2,
$$
and $\Gamma$ is the subgroup generated by $(a_k,b_k)$, $k=1,...,4$
(see, for example, \cite{Borc}). Considering $M$ as the quotient
${\Bbb C}^2/{\Gamma}$, every left-invariant object on ${\Bbb C}^2$
descends to a globally defined object on $M$.

We identify ${\Bbb C}^2$ with ${\Bbb R}^4$ by
$(z=x+iy,w=u+iv)\to(x,y,u,v)$ and set
$$
A_1=\frac{\partial}{\partial x}-x\frac{\partial}{\partial
u}+y\frac{\partial}{\partial v},\quad A_2=\frac{\partial}{\partial
y}-y\frac{\partial}{\partial u}-x\frac{\partial}{\partial v},\quad
A_3=\frac{\partial}{\partial u},\quad A_4=\frac{\partial}{\partial
v}.
$$
These form a basis for the space of left-invariant vector fields on
${\Bbb C}^2$. We note that their Lie brackets are
$$
[A_1,A_2]=-2A_4,\quad [A_i,A_j]=0 %~~\rm{otherwise}
$$
for all other $i,j$. It follows that the group ${\Bbb C}^2$ defined
above is solvable.

 Denote by $g$ the left-invariant Riemannian metric on $M$ for which the basis
$A_1,...,A_4$ is orthonormal.

By \cite{Has}, every complex structure on $M$ is induced by a
left-invariant complex structure on ${\Bbb C}^2$. It is easy to see
that every such a structure is given by (\cite{D14, M})
$$
JA_1=\varepsilon_1 A_2,\quad JA_3=\varepsilon_2 A_4,\quad
\varepsilon_1,\varepsilon_2=\pm 1.
$$
Denote the complex structure defined by these identity by
$J_{\varepsilon_1,\varepsilon_2}$.

The non-zero covariant derivatives $\nabla_{A_i}A_j$ are
(\cite{DHM})
$$
\nabla_{A_1}A_2=-\nabla_{A_2}A_1=-A_4,\quad
\nabla_{A_1}A_4=\nabla_{A_4}A_1=A_2,\quad
\nabla_{A_2}A_4=\nabla_{A_4}A_2=-A_1.
$$
This implies that the Lie form is
$$
\theta(X)=-2\varepsilon_1g(X,A_3).
$$
Therefore
$$
B=-2\varepsilon_1A_3,\quad \nabla\theta=0.
$$

Denote  the dual basis of $A_1,...,A_4$ by $\alpha_1,...,\alpha_4$.
Thus $\theta=-2\varepsilon_1\alpha_3$.  Any left-invariant $1$-form
$\tau$ is of the form $ \tau=a_1\alpha_1+...+a_4\alpha_4$, where
$a_1,...,a_4$ are real constants. We have
$d\alpha_1=d\alpha_2=d\alpha_3=0$, $d\alpha_4=2\alpha_1\wedge
\alpha_2$.  The $3$-form ${\mathcal T}$ corresponding to $\tau$ is $
{\mathcal
T}=a_4\alpha_{123}-a_3\alpha_{124}+a_2\alpha_{134}-a_1\alpha_{234}$
where $\alpha_{ijk}=\alpha_i\wedge\alpha_j\wedge\alpha_k.$ Hence the
form $d\theta-d\tau-\imath_{B}{\mathcal
T}=2(\varepsilon_1-1)a_4\alpha_1\wedge\alpha_2-2\varepsilon_1(a_2\alpha_1\wedge\alpha_4-a_1\alpha_2\wedge\alpha_4)$
is of type $(1,1)$ w.r.t. $J$ if and only if $a_1=a_2=0$. Let $D$ be
the metric connection with skew-symmetric torsion $T$ determined by
the form $\tau$. Then the non-zero covariant derivatives
$D_{A_i}A_j$ are given in the following table
$$
\begin{array}{c}
D_{A_1}A_2=\frac{1}{2}a_4A_3-(1+\frac{1}{2}a_3)A_4,\quad
D_{A_1}A_3=-\frac{1}{2}a_4A_2,\quad
D_{A_1}A_4=(1+\frac{1}{2}a_3)A_2\\[6pt]

D_{A_2}A_1=-\frac{1}{2}a_4A_3+(1+\frac{1}{2}a_3)A_4,\quad
D_{A_2}A_3=\frac{1}{2}a_4A_1,\quad
D_{A_2}A_4=-(1+\frac{1}{2}a_3)A_1\\[6pt]

D_{A_3}A_1=\frac{1}{2}a_4A_2,\quad
D_{A_3}A_2=-\frac{1}{2}a_4A_1,\quad

D_{A_4}A_1=(1-\frac{1}{2}a_3)A_2,\quad
D_{A_4}A_2=-(1-\frac{1}{2}a_3)A_1
\end{array}
$$
Using this table, one can compute the values
${\rho_D}_{ij}={\rho_D}(A_i,A_j)$ of the Ricci tensor. The non-zero
ones are
$$
\begin{array}{c}

\displaystyle{{\rho_D}_{11}={\rho_D}_{22}=-\frac{a_3^2+a_4^2+4}{2}, \quad {\rho_D}_{33}=-\frac{a_4^2}{2}, \quad {\rho_D}_{34}=a_4(1+\frac{1}{2}a_3)}, \\[6pt]

\displaystyle{{\rho_D}_{43}=-a_4(1-\frac{1}{2}a_3), \quad
{\rho_D}_{44}=2-\frac{1}{2}a_3^2}.

\end{array}
$$
Also, the non-zero values of
${\rho^{\ast}_D}_{ij}={\rho^{\ast}_D}(A_i,A_j)$ of the $\ast$-Ricci
tensor are
$$
\begin{array}{c}

\displaystyle{{\rho^{\ast}_D}_{11}=-{\rho^{\ast}_D}_{22}=\frac{a_3^2+a_4^2+12}{4},
\quad
{\rho^{\ast}_D}_{12}=-{\rho^{\ast}_D}_{21}=-\varepsilon_1\varepsilon_2a_4}.

\end{array}
$$
It follows that the identity
$$
\rho_D(Z,B)-\rho^{\ast}_D(Z,B)-\rho_{D}(Z,\tau)+\rho_{D}^{\ast}(Z,\tau)=0.
$$
 is satisfied for every $Z\in TM$ if and only if
$$
(1-\varepsilon_1)a_4=0,\quad (1-\varepsilon_1)(2-a_3)a_4=0.
$$
Clearly, if $\varepsilon_1=1$, the system is satisfied for any $a_3$
and $a_4$, and if $\varepsilon_1=-1$, the solution of the system is
$a_4=0$, $a_3$-arbitrary. Thus, on a Kodaira surface, there are many
metric connections with totally skew-symmetric torsion for which the
complex structures $J_{\varepsilon_1,\varepsilon_2}$ are harmonic
maps from $(M,g)$ into its twistor space $({\mathcal Z}, h_t^D)$ .

\smallskip

The next statement is a weaker version of \cite[Proposition
1.3]{Pon}.

\begin{cor}\label{unique}
Let $J_1$ and $J_2$ be two  Hermitian structures on a (connected)
Riemannian four-manifold $(M,g)$. Suppose that $J_1$ and $J_2$ have
the same Lee form. If $J_1$ and $J_2$ coincide on an open subset or,
more-generally, if they coincide to infinite order at a point, they
coincide on the whole manifold $M$.
\end{cor}

\noindent {\bf Proof}. The complex structures $J_1$ and $J_2$
determine the same orientation on $M$, hence the same positive
twistor space ${\mathcal Z}$. Also, $J_1$ and $J_2$ determine the
same Bismut-Strominger connection $D$. If $h_t$ is the metric on
${\mathcal Z}$ defined by means of $D$, the maps ${\mathfrak J}_1$
and ${\mathfrak J_2}$ of $(M,g)$ into $({\mathcal Z},h_t)$ are
harmonic. Thus, the result follows from the uniqueness theorem for
harmonic maps \cite{S}.

\smallskip

\noindent {\bf Remark}. In fact, M. Pontecorvo \cite{Pon} has proved
that the conclusion of the above corollary holds without the
assumption on the Lee forms. The idea of his proof is inspired by
the theory of pseudo-holomorphic curves.

\medskip

\noindent {\bf Acknowledgement}. I would like to thank Christian
Yankow who checked the curvature computation in Example 1 by a
computer programme and Kamran Shakoor who checked Example 2.

\end{document}